\magnification=1200
\def\ine{{1\over 2\pi} \int_0^{2\pi}}
\def\ha{{1\over 2}}
\def\d{{\rm d}}
\def\L{{\cal L}}
\hfill{Preprint SISSA 25/98/FM}
\vskip 4cm
\centerline{\bf FLAT PENCILS OF METRICS AND FROBENIUS MANIFOLDS}
\medskip
\centerline{\bf Boris Dubrovin}
\medskip
\centerline{SISSA, Trieste, Italy}
\bigskip
{\bf Abstracts} This paper is based on the author's talk at 1997 Taniguchi
Symposium ``Integrable Systems and Algebraic Geometry''. We consider an
approach to the theory of Frobenius manifolds  based on the geometry of
flat pencils of contravariant metrics. It is shown that, under certain
homogeneity assumptions, these two objects are identical. The flat pencils
of contravariant metrics on a manifold $M$ appear naturally in the
classification of bihamiltonian structures of hydrodynamics type on the
loop space $L(M)$. This elucidates the relations between Frobenius
manifolds and integrable hierarchies.
\vfill\eject

{\bf Introduction}
\medskip
Let $M$ be $n$-dimensional smooth manifold.
\smallskip
{\bf Definition 0.1.} A symmetric bilinear form $(~,~)$ on $T^*M$ is
called {\it contravariant metric} if it is invertible on an open dense
subset $M_0\subset M$.
\medskip
In local coordinates $x^1$, \dots, $x^n$ a contravariant metric is
specified by the components (a $(2,0)$-tensor)
$$
\left( dx^i, dx^j\right) = g^{ij}(x), ~i,\, j =1, \dots, n.
\eqno(0.1)
$$
On $M_0$ the inverse matrix $\left( g_{ij}(x)\right) = \left(
g^{ij}(x)\right)^{-1}$ determines a metric in the usual sense
$$
ds^2 = g_{ij}(x) dx^i dx^j
\eqno(0.2)
$$
(not necessarily positive definite). Here and below a summation over
repeated indices is assumed.
\smallskip
{\bf Definition 0.2.} The {\it contravariant Levi-Civita connection}
for the metric $(~,~)$ is determined by a collection of $n^3$ functions
$\Gamma_k^{ij}(x)$ defined for any coordinate patch on $M$ such that on
$M_0$ 
$$
\Gamma_k^{ij}(x) = - g^{is}(x)\Gamma_{sk}^j(x)
\eqno(0.3)
$$
where $\Gamma_{sk}^j(x)$ is the Levi-Civita connection for the metric
(0.2).
\smallskip
{\bf Lemma 0.1.} {\it The coefficients $\Gamma_k^{ij}(x)$ of the
contravariant Levi-Civita connection are determined uniquely on $M_0$
from the system of linear equations
$$
\eqalignno{
g^{is}\Gamma_s^{jk} & = g^{js}\Gamma_s^{ik} & 
(0.4)
\cr
\Gamma_k^{ij} + \Gamma_k^{ji} &= \partial_k g^{ij}. &
(0.5)
\cr}
$$
}

Here 
$$
\partial_k ={\partial\over \partial x^k}.
$$

Proof. On $M_0$ the symmetry condition $\Gamma_{ij}^k =\Gamma_{ji}^k$ of
the Levi-Civita connection reads in the form (0.4). The equation
$\nabla (~,~)=0$ coincides with (0.5). Proof now
 follows from the wellknown theorem of existence and uniqueness
of the Levi-Civita connection. 
\medskip
We emphasize, however, the assumption
of the contravariant connection to be defined on all $M$ but not only on
$M_0$ (the linear system (0.4), (0.5) for $\Gamma^{ij}_k$ degenerates
on $M\setminus M_0$).

Having a contravariant connection one can define the operators $\nabla^u$
of covariant derivatives along any 1-form $u\in T^*M$. For example, for a
covector field $v=v_i dx^i$ we obtain
$$
\nabla^u v = \left( g^{ij} u_i \partial_j v_k + \Gamma_k^{ij}u_i
v_j\right) dx^k.
\eqno(0.6)
$$
Particularly, for $u=dx^i$ the operators $\nabla^i :=\nabla^u$ are related
with the usual covariant derivatives defined on $M_0$ by raising of the
index
$$
\nabla^i =g^{is}\nabla_s.
\eqno(0.7)
$$
If the covector $v_j$ is defined globally on $M$ then $\nabla^i v_j$
will be globally defined (1,1)-tensor on $M$. From this it easily follows
the correctness of the definition of contravariant connection.

\medskip
{\bf Definition 0.3.} A function $f(x)$ is called {\it flat coordinate} of
$(~,~)$ if the differential $\xi = df$ is covariantly constant w.r.t. the
Levi-Civita connection
$$
g^{is} \partial_s \xi_j +\Gamma_j^{is} \xi_s =0, ~i, \, j =1, \dots, n,
\eqno(0.8)
$$
$$
\xi_j =\partial_j f.
$$
\smallskip
{\bf Definition 0.4.} A contravariant metric is said to be {\it flat} 
{\it iff} on $M_0$ there locally exist $n$ independent flat coordinates.
\medskip
Choosing a system of flat coordinates one reduces the matrix (0.1) of the
metric to a constant form and the coefficients $\Gamma_k^{ij}$ of the
Levi-Civita connection to zero.
\smallskip
{\bf Lemma 0.2.} {\it The contravariant metric is flat {\rm iff} the
Riemann curvature tensor
$$
R_l^{ijk} := g^{is} \left( \partial_s \Gamma_l^{jk}
-\partial_l\Gamma_s^{jk}\right) +
\Gamma_s^{ij}\Gamma_l^{sk}-\Gamma_s^{ik}\Gamma_l^{sj}
\eqno(0.9)
$$
identically vanishes.
}

This is a standard fact of differential geometry (see, e.g., [DFN]). It is
important that our formula for the curvature involves only contravariant
components of the metric and of the connection.

We give now our main
\smallskip
{\bf Definition 0.5.} Two contravariant metrics $(~,~)_1$ and $(~,~)_2$
form a {\it flat pencil} if:

1) The linear combination
$$
(~,~)_1 -\lambda (~,~)_2
\eqno(0.10)
$$
for any $\lambda$ is a contravariant metric on $M$.

2) If ${\Gamma_1}_k^{ij}$ and ${\Gamma_2}_k^{ij}$ are the contravariant
Levi-Civita connections for these two metrics then for any $\lambda$ the
linear combination
$$
{\Gamma_1}_k^{ij}-\lambda {\Gamma_2}_k^{ij}
$$
is the contravariant Levi-Civita connection for the metric (0.10).

3). The metric (0.10) is flat for any $\lambda$.

We say that the flat pencil of metrics is {\it quasihomogeneous of the
degree $d$} if there exists a function $\tau$ on $M$ such that the vector
fields
$$
E:=\nabla_1 \tau, ~~E^i =g_1^{is}\partial_s\tau
\eqno(0.11a)
$$
$$
e:=\nabla_2 \tau, ~~e^i = g_2^{is}\partial_s\tau
\eqno(0.11b)
$$
satisfy the following properties
$$
[e,E]=e
\eqno(0.12)
$$
$$
\L_E (~,~)_1 =(d-1) (~,~)_1
\eqno(0.13)
$$
$$
\L_e (~,~)_1 = (~,~)_2
\eqno(0.14)
$$
$$
\L_e(~,~)_2 =0.
\eqno(0.15)
$$
\medskip
\def\deg{\mathop{\rm deg}}
\def\C{{\bf C}}
\smallskip
{\bf Definition 0.6.} A {\it Frobenius algebra} is a pair $\left( A,
~<~,~>\right)$ where $A$ is a commutative associative algebra (over
{\bf R} or {\bf
C}) with a unity and $<~,~>$ stands for a symmetric
non-degenerate {\it invariant} bilinear form on $A$. The invariance
means validity of the following identity
$$<a\, b, c> = < a, b\, c>
\eqno(0.16)
$$
for arbitrary 3 vectors $a, ~b, ~c \in A$.
\smallskip
{\bf Definition 0.7.} The Frobenius algebra is called {\it graded} if a
linear operator $Q: A\to A$ and a number $d$ are defined such that
$$
Q(a \, b) = Q(a) b + a Q(b),
\eqno(0.17a)
$$
$$
<Q(a),b> + <a, Q(b)> = d <a,b>
\eqno(0.17b)
$$ 
for any $a, ~b \in A$. The operator $Q$ is called {\it grading operator}
and the number $d$ is called {\it charge} of the Frobenius algebra. 
In 
the case of diagonalizable grading operators we
may
assign degrees to the eigenvectors $e_\alpha$ of $Q$
$$
\deg (e_\alpha) =q_\alpha ~~{\rm if} ~~ Q(e_\alpha) =q_\alpha e_\alpha.
\eqno(0.18)
$$
Then the usual property of the degree of the product of homogeneous
elements of the algebra holds true
$$
\deg\,(ab) =\deg\,a \deg\,b.
$$
Besides, $<a,b>$ can be nonzero only if $\deg\,a + \deg\, b=d$
where $d$ is the charge.

We will consider also graded Frobenius algebras
$\left( A, <~,~>\right)$ over
graded commutative associative rings $R$. In this case we have two grading
operators $Q_R: R\to R$ and $Q_A:A\to A$
satisfying the properties
$$
\eqalignno{Q_R(\alpha\beta) & = Q_R (\alpha)\beta + \alpha Q_R(\beta),~
\alpha, ~\beta \in R & (0.19a)
\cr
Q_A(ab) & = Q_A(a)b+aQ_A(b), ~a, ~b \in A & (0.19b)
\cr
Q_A(\alpha a) &= Q_R(\alpha)a + \alpha Q_A(a), ~\alpha\in R, ~a\in A
& (0.19c)
\cr
Q_R<a,b> +d\, <a,b> &=<Q_A(a),b> + <a, Q_A(b)>, ~a, ~b\in A.
& (0.19d)
\cr}
$$
As above the number $d$ is called {\it the charge} of the graded Frobenius
algebra over the graded ring.
\medskip
{\bf Definition 0.8.} (Smooth, analytic) {\it Frobenius structure} 
on the manifold $M$ is a structure of Frobenius algebra on the tangent spaces 
$T_tM=\left(A_t, <~,~>_t\right)$
depending (smoothly, analytically) on the point $t$. This structure must 
satisfy the following axioms.
\smallskip
\noindent{\bf FM1}. The metric on $M$ induced by the invarint bilinear form
$<~,~>_t$ is flat. Denote $\nabla$ is the Levi-Civita connection for the metric
$<~,~>_t$. The unity vector field $e$ must be covariantly constant,
$$\nabla \, e=0.
\eqno(0.20)
$$
As above we use here the word `metric' as a synonim of a symmetric
nondegenerate
bilinear form on $TM$, not necessarily of a positive one. Flatness of the
metric,
i.e., vanishing of the Riemann curvature tensor, means that locally a
{\it system of flat coordinates} $(t^1, \dots, t^n)$ exists such that the 
matrix $<\partial_\alpha , \partial_\beta>$
of the metric in these coordinates becomes constant.
\medskip
\noindent{\bf FM2}. Let $c$ be the following symmetric trilinear form on $TM$
$$c(u,v,w) := <u\cdot v, w>.
\eqno(0.21)
$$
The four-linear form
$$
(\nabla_z c)(u,v,w), ~u,v,w,z\in TM
\eqno(0.22)
$$
must be also symmetric.
\medskip
Before formulating the last axiom we observe that the space $Vect(M)$ of 
vector fields on $M$ acquires a structure of a Frobenius algebra over the 
algebra
$Func(M)$ of (smooth, analytic) functions  on $M$. 
\smallskip
\noindent{\bf FM3}. A linear {\it Euler vector field} $E\in Vect(M)$ must
be fixed on $M$,
i.e.,
$$
\nabla\nabla \, E =0.
\eqno(0.23)
$$
The operators
$$
\eqalign{
Q_{Func(M)}&:=E\cr
Q_{Vect(M)}&:={\rm id} + {\rm ad}_E\cr}
\eqno(0.24)
$$
introduce in $Vect(M)$ a structure of graded Frobenius algebra  of a
given charge $d$ over the graded
ring $Func(M)$. 
\medskip
We will now spell out the requirements of Definition 0.8 in the flat
coordinates $t^1$, \dots, $t^n$ of the metric $<~,~>$. Denote
$$
\eta_{\alpha\beta} :=<\partial_\alpha,\partial_\beta>
\eqno(0.25)
$$
(a constant symmetric nondegenerate matrix),
$$
\partial_\alpha\cdot\partial_\beta |_{t} = c_{\alpha\beta}^\gamma(t)
\partial_\gamma.
\eqno(0.26)
$$
Then the components of the trilinear form (0.21) can be locally
represented as the triple derivatives of a function $F(t)$
$$
c\left( \partial_\alpha, \partial_\beta, \partial_\gamma \right) =
 \partial_\alpha \partial_\beta \partial_\gamma F(t)
\eqno(0.27).
$$
Associativity of the multiplication in the Frobenius algebra implies
the following {\it WDVV associativity equations} for the function $F(t)$
$$
 \partial_\alpha \partial_\beta \partial_\lambda F(t) \eta^{\lambda\mu}
 \partial_\mu \partial_\gamma \partial_\delta F(t)
=
\partial_\delta \partial_\beta \partial_\lambda F(t) \eta^{\lambda\mu}
 \partial_\mu \partial_\gamma \partial_\alpha F(t), ~\alpha, \, \beta,\,
\gamma,\, \delta =1, \dots, n.
\eqno(0.28)
$$
The axiom {\bf FM3} in these coordinates can be recasted into the
following equivalent form
$$
\L_E c_{\alpha\beta}^\gamma =  c_{\alpha\beta}^\gamma
\eqno(0.29)
$$
$$
\L_E \eta_{\alpha\beta}=(2-d)\eta_{\alpha\beta}.
\eqno(0.30)
$$
From (0.29) one obtains the following quasihomogeneity equation for the
function  $F(t)$
$$
\L_EF(t) =(3-d) F(t) +\ha A_{\alpha\beta}t^\alpha t^\beta + B_\alpha
t^\alpha +C
\eqno(0.31)
$$
with some constants $A_{\alpha\beta}$, $B_\alpha$, $C$.

Finally, {\bf FM1} means that the unity vector field is constant in the
coordinates $t$. Usually the flat coordinates are  chosen in such a way
that
$$
e =\partial / \partial t^1.
\eqno(0.32)
$$
The equation (0.23) means that the matrix $\nabla E$ is constant in the
flat coordinates $t$.

The main aim of the present paper is to prove that any Frobenius manifold
carries a natural quasihomogeneous linear pencil of metrics and, under
certain nondegeneracy assumption, to prove also the converse statement.
\bigskip
\centerline{\bf 1. From Frobenius manifolds to flat pencils}
\medskip
We put $(~,~)_2=<~,~>$ (as a bilinear form on the cotangent bundle) and we
define a new bilinear form on the cotangent bundle
$$
(\omega_1, \omega_2)_1 \equiv (\omega_1, \omega_2) = i_E (\omega_1 \cdot
\omega_2).
\eqno(1.1)
$$
Here $i_E$ is the operator of inner product (i.e., the contraction of the
vector field $E$ with a 1-form).This metric was found in [Du2] to give
bihamiltonian structure of the integrable hierarchies of [Du1] describing
coupling of a given matter sector of a 2D topological field theory to
topological gravity. It was called in [Du4] {\it
intersection form} of the Frobenius manifold.
 \smallskip
{\bf Theorem 1.1.} {\it The metrics $(~,~)$ and $<~,~>$ on a Frobenius
manifold form a flat pencil quasihomogeneous of the degree $d$.}

Proof essentially follows [Du4] (we only relax the assumptions about the
Frobenius manifold not requiring diagonalizability of the tensor $\nabla
E$).

In the flat coordinates $t^1$, \dots, $t^n$ for $<~,~>$ the components of
the bilinear form (1.1) are given by the formula
$$
\eqalign{
g^{\alpha\beta}(t) &=\left(dt^\alpha, dt^\beta\right) =E^\epsilon(t)
c_\epsilon^{\alpha\beta}(t)
\cr
&= R_\epsilon^\alpha F^{\epsilon\beta}(t) +F^{\alpha\epsilon}(t)
R_\epsilon^\beta +A^{\alpha\beta}\cr}
\eqno(1.2)
$$
where
$$
R_\beta^\alpha ={d-1\over 2}\delta_\beta^\alpha +\left( \nabla
E\right)_\beta^\alpha
\eqno(1.3)
$$
$$
c_\gamma^{\alpha\beta}(t) =\eta^{\alpha\epsilon} c_{\epsilon\gamma}^\beta
(t)
$$
$$
F^{\alpha\beta}(t)
=\eta^{\alpha\lambda}\eta^{\beta\mu}\partial_\lambda\partial_\mu F(t)
\eqno(1.4)
$$
$$
A^{\alpha\beta} =\eta^{\alpha\lambda}\eta^{\beta\mu}A_{\lambda\mu}
$$
and the matrix $A_{\alpha\beta}$ was defined in (0.31).

From the first line in (1.2) it follows that
$$
g^{\alpha\beta}(t) = \eta^{\alpha\beta} t^1 +\tilde g^{\alpha\beta}(t^2,
\dots, t^n).
$$
So
$$
g^{\alpha\beta}(t)-\lambda \eta^{\alpha\beta} = \eta^{\alpha\beta}
(t^1-\lambda) +\tilde g^{\alpha\beta}(t^2,  
\dots, t^n)
$$
for any $\lambda$ does not degenerate on an open dense subset in $M$.
\smallskip
{\bf Lemma 1.1.} {\it The contravariant Levi-Civita connection for the
metric
$(~,~)-\lambda <~,~>$ is given by the formula
$$
\Gamma_\gamma^{\alpha\beta} (t) = c^{\alpha\epsilon}_\gamma (t)
R_\epsilon^\beta.
\eqno(1.5)
$$
}

Proof. Differentiating (1.2) w.r.t. $t^\gamma$ we obtain
$$
\partial_\gamma g^{\alpha\beta} = R_\epsilon^\alpha
c^{\epsilon\beta}_\gamma + c_\gamma^{\alpha\epsilon}
R_\epsilon^\beta=\Gamma_\gamma^{\beta\alpha}+\Gamma_\gamma^{\alpha\beta}.
$$
This proves (0.5).
The second equation (0.4) follows immediately from associativity of
the multiplication on $T^*M$. Lemma is proved.
\medskip
To finish the proof of Theorem it remains to show that the curvature of
the pencil of the metrics vanishes identically in $\lambda$. First observe
that the terms with derivatives of $\Gamma$ in
$$
\left( g^{\alpha\epsilon} - \lambda \eta^{\alpha\epsilon}\right)
\left( \partial_\epsilon \Gamma_\delta^{\beta\gamma} -
\partial_\delta \Gamma_\epsilon^{\beta\gamma} \right)
+ \Gamma_\epsilon^{\alpha\beta} \Gamma_\delta^{\epsilon\gamma}
-\Gamma_\epsilon^{\alpha\gamma} \Gamma_\delta^{\epsilon\beta}
$$
vanish
due to equality of mixed derivatives
$$
 \partial_\epsilon c_\delta^{\beta\gamma} =
\partial_\delta c_\epsilon^{\beta\gamma}.
$$
The remaining terms vanish due to associativity. 

Finally, we put
$$
\tau = \eta_{1\alpha} t^\alpha
\eqno(1.6)
$$
assuming that the coordinate $t^1$ is chosen as in (0.32). Then (0.11)
immediately follows. The equations (0.12) -- (0.15) follow from (0.29),
(0.30) and from
$$
\L_E e =-e
\eqno(1.7)
$$
(a consequence of (0.29)).
Theorem is proved.

\bigskip
\centerline{\bf 2. From flat pencils to Frobenius manifolds}
\medskip
We begin with simple
\smallskip
{\bf Lemma 2.1.} {\it The functions
$$
\Delta^{ijk}(x) = g_2^{js} {\Gamma_1}_s^{ik}-g_1^{is}{\Gamma_2}_s^{jk}
\eqno(2.1)
$$
are components of a rank three tensor (i.e., of a trilinear form on
$T^*M$). Two flat metrics $g_1^{ij}$ and $g_2^{ij}$ can be simultaneously
reduced to constant form {\rm iff} $\Delta^{ijk}=0$.
}

Proof. It is wellknown that the difference of usual Christoffel
coefficients of two affine connections 
$$
{\Gamma_1}_{st}^k -{\Gamma_2}_{st}^k
\eqno(2.2)
$$
is a tensor. Contracting this with $g_1^{is}g_2^{jt}$ we obtain the tensor
(2.1). Two metrics are simultaneously reducible to a constant form {\it
iff} the difference (2.2) vanishes. Lemma is proved.
\medskip
We will also consider a (2,1)-tensor
$$
\Delta_i^{jk} ={g_2}_{is} \Delta^{sjk}
\eqno(2.3)
$$
defined on the open subset $M_0\subset M$ where the contravariant metric
$(~,~)_2$ does not degenerate. The tensor (2.3) defines a bilinear
operation
$$
T^*M_0 \times T^*M_0 \to T^*M_0
\eqno(2.4a)
$$
$$
\eqalign{
&(u,v)\mapsto \Delta(u,v), 
\cr
&\Delta(u,v)_k =u_i v_j \Delta_k^{ij}(x) ~{\rm for} ~u, \, v \in T_x^*M_0.
\cr}
\eqno(2.4b)
$$
\smallskip
{\bf Lemma 2.2.} {\it For a flat pencil of metrics the tensor (2.3)
satisfies the following properties
$$
\left( \Delta(u,v),w\right)_1 = \left( u, \Delta(w,v)\right)_1
\eqno(2.5)
$$
$$
\left( \Delta(u,v),w\right)_2 = \left( u, \Delta(w,v)\right)_2
\eqno(2.6)
$$
$$
\Delta(\Delta(u,v),w) = \Delta(\Delta(u,w),v)
\eqno(2.7)
$$
$$
\nabla_2^u \Delta(v,w) - \nabla_2^v \Delta(u,w) =\Delta\left( \nabla_2^u v
-
\nabla_2^v u, w\right).
\eqno(2.8)
$$
Here $\nabla_2^u$ are the covariant derivative operators (0.6) for the
second metric $(~,~)_2$, $u$, $v$, $w$ are arbitrary 1-forms on $M_0$.
}

Quasihomogeneity of the flat pencil is not assumed. Note that, due to
(2.7), the algebra structure on  $T^*M_0$ will not be associative but {\it
right-symmetric}.

Proof (see [Du4], Appendix D). Let us choose a system of flat coordinates
$x^1$, \dots, $x^n$ for the metric $(~,~)_2$. In these coordinates we have
$$
{\Gamma_2}_k^{ij}=0, ~{\Gamma_1}_k^{ij} = \Delta_k^{ij}.
$$
From the definition it follows that $\Delta_k^{ij}$ will also coincide
with the contravariant Levi-Civita connection for all the linear pencil
$(~,~)_1-\lambda (~,~)_2$ with an arbitrary $\lambda$. Writing the
symmetry condition (0.4)
$$
\left( g^{is}_1 -\lambda g_2^{is}\right) \Delta_s^{jk}=
\left( g^{js}_1 -\lambda g_2^{js}\right) \Delta_s^{ik}
$$
we obtain (2.5) and (2.6). Vanishing of the curvature (0.9) of the pencil
gives the equations
$$
\left( g^{is}_1 -\lambda g_2^{is}\right) \left(\partial_s\Delta_l^{jk}
-\partial_l \Delta_s^{jk}\right) +\Delta_s^{ij} \Delta_l^{sk}
-\Delta_s^{ik}\Delta_l^{sj}=0 ~{\rm for~any}~i,\,j,\,k,\,l.
$$
Vanishing of the linear in $\lambda$ term implies
$$
\partial_s\Delta_l^{jk} 
-\partial_l \Delta_s^{jk}=0.
$$
This coincides with (2.8). Vanishing of the remaining terms gives (2.7).
Lemma is proved.
\medskip
{\bf Lemma 2.3.} {\it For a quasihomogeneous flat pencil the following
equations hold true
$$
\nabla_2\nabla_2\tau=0
\eqno(2.9)
$$
$$
\nabla_2\nabla_2 E=0.
\eqno(2.10)
$$
}

Proof. We have
$$
0=\L_e {g_1}_{ij} =2 {\nabla_2}_i {\nabla_2}_j \tau.
$$
This proves (2.9). From (0.12) and (0.13) it follows
$$
\L_E (~,~)_2 =(d-2) (~,~)_2.
\eqno(2.11)
$$
So the vector field $E$ generates the one-parameter group of linear
conformal transformations of the metric $(~,~)_2$. This proves (2.10).
\medskip
{\bf Corollary 2.1.} {\it The eigenvalues of the matrix
$$
{\nabla_2}_i E^j(x)
\eqno(2.12)
$$
do not depend on the point of the manifold.}
\smallskip
{\bf Definition 2.1.} A  quasihomogeneous flat pencil is said to be {\it
regular} if the (1,1)-tensor
$$
R_i^j = {d-1\over 2}\delta_i^j + {\nabla_2}_i E^j
\eqno(2.13)
$$
does not degenerate on $M$.
\smallskip
{\bf Theorem 2.1.} {\it Let $M$ be a manifold carrying a regular
quasihomogeneous flat pencil. Denote $M_0\subset M$ the
subset of $M$ where  the metric $(~,~)_2$ is invertible.
Define the
multiplication
of 1-forms on $M_0$ 
putting
$$
u\cdot v
:= \Delta\left( u, R^{-1}v\right).
\eqno(2.14)
$$
Then there exists a unique Frobenius structure on $M$ such that 
$$
<~,~> = (~,~)_2,
\eqno(2.15)
$$
the multiplication of tangent vectors is $<~,~>$-dual to the product
(2.14), the unity and the Euler vector fields have the form (0.11a) and
(0.11b)
resp., and the intersection form is equal to $(~,~)_1$.
}

Proof. Let us choose flat coordinates $t^1$, \dots, $t^n$ for the metric
$(~,~)_2$. The components of the metric in these coordinates are given by
a constant symmetric invertible matrix
$$
\eta^{\alpha\beta} : =\left( dt^\alpha, dt^\beta\right)_2.
\eqno(2.16)
$$
We also denote
$$
g^{\alpha\beta}(t) : =\left( dt^\alpha, dt^\beta\right)_1
$$
and
$$
K_\alpha^\beta := \partial_\alpha E^\beta.
\eqno(2.17)
$$
This matrix is constant due to Lemma 2.3. The components of
the contravariant Levi-Civita
connection for the metric $g^{\alpha\beta}$ in these coordinates we denote
$\Gamma_\gamma^{\alpha\beta}$. Recall that in this coordinate system
$$
\Delta_\gamma^{\alpha\beta} =\Gamma_\gamma^{\alpha\beta}.
\eqno(2.18)
$$ 

\smallskip
{\bf Lemma 2.4.} {\it The vector field $e$ is constant in the coordinates
$t^\alpha$. It is an eigenvector of the operator (2.17) with the
eigenvalue $1-d$.
}

Proof. Constancy of $e$ follows from (2.9). Let us normalize the choice of
the flat coordinates requiring that
$$
e^\alpha =\eta^{\alpha\,n}.
\eqno(2.19)
$$
In these coordinates
$$
\tau= t^n + {\rm const}.
\eqno(2.20)
$$
So
$$
E^\alpha(t) = g^{\alpha\, n}(t), ~\alpha=1, \dots, n.
\eqno(2.21)
$$
From (0.12) we obtain
$$
\eta^{n\,\epsilon} K_\epsilon^\alpha = \eta^{\alpha\, n}.
$$
Using (2.14) we obtain
$$
\eta^{n\,\epsilon} K_\epsilon^\alpha + \eta^{\alpha\epsilon}
K_{\epsilon}^n = (2-d) \eta^{n\, \alpha}.
$$
Hence
$$
\eta^{\alpha\epsilon} K_\epsilon^n =(1-d) \eta^{n\, \alpha}.
$$
Lowering the index $\alpha$ we prove Lemma.
\medskip
We will use also below the choice (2.19) of the flat coordinate $t^n$.
Then
$$
K^n_\alpha =(1-d) \delta_\alpha^n.
\eqno(2.22)
$$
\smallskip
{\bf Lemma 2.5.} {\it In the coordinates $t^\alpha$
$$
\eqalignno{
\Delta_\beta^{\alpha\, n} &={1-d\over 2} \delta^\alpha_\beta & (2.23)
\cr
\Delta_\beta^{n\,\alpha} & ={d-1\over 2}\delta^\alpha_\beta +
K_\beta^\alpha & (2.24)
\cr}
$$
}

Proof. From (0.13) it follows that
$$
\L_E g_{\alpha\beta} =(1-d) g_{\alpha\beta}.
$$
Using Christoffel formula one obtains
$$
\eqalign{
\Gamma_{\alpha\beta}^n &={1\over 2} g^{n\,\epsilon} \left( \partial_\alpha
g_{\epsilon\beta} + \partial_\beta g_{\alpha\epsilon} - \partial_\epsilon
g_{\alpha\beta}\right)
\cr
&= -{1\over 2} \L_E g_{\alpha\beta} ={d-1\over 2} g_{\alpha\beta}.
\cr}
$$
Raising the index we obtain
$$
\Gamma_\beta^{\alpha\,n} ={1-d\over 2} \delta_\beta^\alpha.
$$ 
Due to (2.18) this proves (2.23). Using the equation (0.5)
$$
\Gamma_\gamma^{n\, \alpha} + \Gamma_\gamma^{\alpha\, n}=\partial_\gamma
g^{\alpha\, n} = K^\alpha_\gamma
$$
we obtain (2.24). Lemma is proved.
\medskip
{\bf Lemma 2.6.} 
$$
\L_E \Delta_\gamma^{\alpha\beta} = (d-1) \Delta_\gamma^{\alpha\beta}
\eqno(2.25)
$$
$$
\L_e \Delta_\gamma^{\alpha\beta} =0.
\eqno(2.26)
$$

Proof. Denote
$$
\eqalign{
\tilde\Gamma_\gamma^{\alpha\beta} &:= \L_E \Delta_\gamma^{\alpha\beta}
+(1-d) \Delta_\gamma^{\alpha\beta}
\cr
&\equiv \partial_E \Gamma_\gamma^{\alpha\beta} -K_\epsilon^\alpha
\Gamma_\gamma^{\epsilon\beta} -\Gamma_\gamma^{\alpha\epsilon}
K_\epsilon^\beta + K_\gamma^\epsilon \Gamma_\epsilon^{\alpha\beta} +
(1-d)\Gamma_\gamma^{\alpha\beta}\cr}
$$
Differentiating the equations
$$
\eqalign{
\Gamma_\gamma^{\alpha\beta} + \Gamma_\gamma^{\beta\alpha} &=
\partial_\gamma g^{\alpha\beta}
\cr
g^{\alpha\epsilon} \Gamma_\epsilon^{\beta\gamma} &= g^{\beta\epsilon}
\Gamma_\epsilon^{\alpha\gamma}
\cr}
\eqno(2.27)
$$
along $E$ we obtain, after simple calculations,
$$
\eqalign{
\tilde\Gamma_\gamma^{\alpha\beta} + \tilde\Gamma_\gamma^{\beta\alpha} &=
0
\cr
g^{\alpha\epsilon} \tilde\Gamma_\epsilon^{\beta\gamma} &=
g^{\beta\epsilon}
\tilde\Gamma_\epsilon^{\alpha\gamma}.
\cr}
\eqno(2.28)
$$
Since the system (2.27) has unique solution for given $g^{\alpha\beta}$,
the correspondent linear homogeneous system (2.28) has only trivial
solution $\tilde\Gamma_\gamma^{\alpha\beta}=0$. This proves (2.25). The
equation (2.26) can be proved in a similar way.
\medskip
{\bf Corollary 2.2.} {\it Let $u$, $v$ be two 1-forms covariantly constant
w.r.t. $\nabla_2$. Then
the multiplication
$$
(u,v)\mapsto \Delta(u,v)
$$
on $T^*M$ satisfies the equations
$$
\Delta(u,v) + \Delta(v,u) =\d (u,v)
\eqno(2.29)
$$
$$
\Delta(R(u), v) + \Delta(u, R(v)) = \d (u, R(v)).
\eqno(2.30)
$$
}

We denote by Roman `d' the differential of a function on $M$ to avoid
confusion with the charge $d$ in the axiom {\bf FM3}.

Proof. The first equation is due to the first
line in (2.27) together with (2.18). The second one is a spelling of
(2.25).
\medskip
Let us now fix a point $t\in M_0$. We denote
$$
V=T_{t_0}^*M.
$$
The linear operator
$$
\Lambda: V\to V, ~~\Lambda = {d-2\over 2}{\bf 1} + K =-{1\over 2} {\bf 1}
+ R
\eqno(2.31)
$$
is skew-symmetric w.r.t. $(~,~)_2 = <~,~>$
$$
<\Lambda u, v> + <u,\Lambda v> =0.
\eqno(2.32)
$$
Let
$$
V=\oplus_\lambda V_\lambda
\eqno(2.33)
$$
be the root decomposition of the space $V$ w.r.t. the root subspaces of
the operator $\Lambda$. The following elementary statement is wellknown
\smallskip
{\bf Lemma 2.7.} {\it The root subspaces $V_\lambda$ and $V_\mu$ are
$<~,~>$-orthogonal if $\lambda+\mu\neq 0$. The pairing
$$
<~,~>: \, V_\lambda \times V_{-\lambda} \to \C
\eqno(2.34)
$$
does not degenerate.
}
\medskip
By the moment we have not used the regularity condition
$$
\det R\neq 0.
\eqno(2.35)
$$
If this condition holds true then 
$$
V_{-{1\over 2}} =0.
\eqno(2.36)
$$
Particularly, this implies that
$$
d\neq 1.
\eqno(2.37)
$$
{\bf Lemma 2.8.} {\it The multiplication (2.14) on $V$ is commutative.}

Proof. From (2.30) we derive the folowing property of the multiplication
$$
u\cdot R(v) + R(u)\cdot v =\d (u,v).
\eqno(2.38)
$$
To take the differential in the r.h.s. of the equation we continue
the covectors $u$, $v$ in a small neighbourhood of the point $t_0$
as $\nabla_2$-constant 1-forms.
It suffices to prove Lemma for $u\in V_\lambda$, $v\in V_\mu$.

Case 1: $\lambda+\mu + 1\neq 0$. Let first $u$ and $v$ be the eigenvectors
of $\Lambda$ with the eigenvalues $\lambda$ and $\mu$ resp. Then
$$
R(u) =\left( \ha + \lambda\right) u, ~~R(v) =\left( \ha + \mu \right) v.
$$
From (2.38) we obtain
$$
(1+\lambda+\mu) u\cdot v = \d (u,v).
\eqno(2.39)
$$
This proves that $u\cdot v=v\cdot u$.
Let $u^{(k)}$, $v^{(l)}$ be the adjoint vectors for the eigenvectors
$u$ and $v$ of the heights $k$, $l$ resp., i.e.,
$$
\Lambda(u^{(k)}) = \lambda u^{(k)} + u^{(k-1)}, ~
\Lambda(v^{(l)}) = \mu v^{(l)} + v^{(l-1)},
$$
$$
u^{(0)}=u, ~v^{(0)}=v, ~u^{(-1)}=v^{(-1)}=0.
$$
We use induction w.r.t. the sum of the heights $k+l$. Substituting in
(2.38)
 $u\mapsto u^{(k)}$, $v\mapsto v^{(l)}$ we obtain
$$
(1+\lambda+\mu) u^{(k)}\cdot v^{(l)} + u^{(k)}\cdot v^{(l-1)} + u^{(k-1)}
\cdot v^{(l)} = \d (u^{(k)}, v^{(l)}).
\eqno(2.40)
$$
By induction 
$$
 u^{(k)}\cdot v^{(l-1)}=  v^{(l-1)} \cdot  u^{(k)}, ~
 u^{(k-1)}
\cdot v^{(l)}=  v^{(l)} \cdot  u^{(k-1)}.
$$
This proves commutativity of $u^{(k)}$ and $v^{(l)}$.

Case 2. $\lambda+\mu+1=0$. Again we use induction w.r.t. the sum of the
weights. Let $u$, $v$ be two eigenvectors. From (2.39) one obtains
$$
\d (u,v) =0.
$$
Using (2.29) we conclude that
$$
\Delta(u,v) + \Delta(v,u) =0.
$$
Hence
$$
u\cdot v ={\Delta(u,v)\over \ha +\lambda} = - {\Delta(v,u)\over \ha
+\lambda} = {\Delta(v,u)\over \ha + \mu} = v\cdot u.
$$
Now we prove commutativity of adjoint vectors. From the definition we
have
$$
\eqalign{
u^{(k)} \cdot v^{(l)} & = {\Delta\left( u^{(k)},v^{(l)}\right)\over \ha
+\mu} -{u^{(k)}\cdot v^{(l-1)}\over \ha +\mu}
\cr
v^{(l)}\cdot u^{(k)} &= {\Delta\left( v^{(l)}, u^{(k)}\right)\over \ha
+\lambda} -{v^{(l)}\cdot u^{(k-1)}\over \ha +\lambda}.
\cr}
$$
So
$$
u^{(k)}\cdot v^{(l)} - v^{(l)} \cdot u^{(k)} =
-{\d\left(u^{(k)},v^{(l)}\right) \over \ha +\lambda}
+
{u^{(k)}\cdot v^{(l-1)} + v^{(l)} \cdot u^{(k-1)}\over \ha +\lambda}.
$$
Applying (2.40) we derive that
$$
\d\left( u^{(k)}, v^{(l)}\right) = u^{(k)} \cdot v^{(l-1)} + v^{(l)} \cdot
u^{(k-1)}.
$$
Lemma is proved.
\medskip
End of the proof of Theorem. We obtained a symmetric multiplication
on the cotangent planes $T_t^*M$
$$
\left( dt^\alpha, dt^\beta\right) \mapsto dt^\alpha\cdot dt^\beta
=:c^{\alpha\beta}_\gamma(t) dt^\gamma
\eqno(2.41)
$$
where the coefficients $c_\gamma^{\alpha\beta}(t)$ are defined by this
equation. The 1-form $dt^n=d\tau$ is the unity of this multiplication.
Indeed, due to (2.23) 
$$
\Delta(dt^\alpha, dt^n) = {1-d\over 2} dt^\alpha.
$$
But the 1-form $dt^n$ is an eigenvector of $\Lambda$ with the eigenvalue
$-d/2$ (this follows from (2.22)). So
$$
dt^\alpha \cdot dt^n = dt^\alpha
$$
for any $\alpha$. Associativity of the multiplication follows from the
right-symmetry property (2.7) and from the commutativity.

By duality we obtain a commutative associative multiplication on $T_tM$
$$
\partial_\alpha\cdot \partial_\beta = c_{\alpha\beta}^\gamma(t)
\partial_\gamma
$$
with
$$
c_{\alpha\beta}^\gamma(t)=\eta_{\alpha\epsilon}c^{\epsilon\gamma}_\beta(t).
$$
The vector $e$ of the form (2.19) will be the unity of this
multiplication.
From commutativity of the multiplication and from (2.6) it follows that
the tensor
$$
<\partial_\alpha\cdot \partial_\beta, \partial_\gamma>
$$
is symmetric w.r.t. $\alpha$, $\beta$, $\gamma$. From this and from (2.8)
it follows
that the gradient
$$
\partial_\delta <\partial_\alpha\cdot \partial_\beta, \partial_\gamma>
$$
is symmetric w.r.t. all the four indices. This proves {\bf FM2}. 

The equation (2.10) implies (0.23).
From the
definition of $K$ it follows that
$$
\L_E K_\alpha^\beta =0.
$$
Hence
$$
\L_E c^{\alpha\beta}_\gamma = (d-1) c^{\alpha\beta}_\gamma.
$$
Lowering the index $\alpha$ we obtain
$$
\L_E c_{\alpha\gamma}^\beta =  c_{\alpha\gamma}^\beta.
$$
This proves (0.29). The equation (0.30) follows from (2.11). 
So, we obtained a Frobenius structure on $M_0$.
Finally,
comparing the equation (2.38) with the second line in the formula (1.2)
for
the entries of the intersection form we conclude that the metric
$g^{\alpha\beta}$ coincides with the intersection form of the Frobenius
manifold. Theorem is proved.
\medskip
{\bf Remark.} In some cases the regularity assumption of non 
degenerateness
of the
operator $R={d-1\over 2}{\bf 1} + \nabla_2 E$ can be relaxed. For example,
for $d=1$ the operator $R$ is always degenerate since
$$
R(d\tau)=0.
$$
However, Theorem 2.1 remains valid under the assumption that the root
subspace $V_{-\ha}$ (see Lemma 2.7 above) is exactly one-dimensional.
Indeed, using the above construction we arrive at  multiplication $u\cdot
v$ defined for an arbitrary 1-form $u$ and for any 1-form $v$ that belongs
to the image of $R$. The only 1-form not belonging to the image is
$d\tau$. However, from (2.24) we obtain
$$
\Delta(d\tau, v) = R(v).
$$
So $d\tau$ is {\it left unity} of the multiplication. Defining $v\cdot
d\tau = v$ we obtain the needed Frobenius structure (cf. [DZ1], proof of
Theorem 2.1).
\medskip
{\bf Example 2.1.} Let $W$ be an irreducible finite Coxeter group acting
in
the Euclidean space ${\bf R}^n$. Denote $(~,~)$ the $W$-invariant
Euclidean inner product on ${\bf R}^n$. According to Arnold [Arn] there
exists a unique contravariant metric $(~,~)_1$ on the orbit
space
$$
M=\C^n /W
$$
such that for any two $W$-invariant polynomials $p(x)$, $q(x)$
$$
(dp,dq)_1 = (dp(x), dq(x)).
\eqno(2.42)
$$
Here we consider $dp$, $dq$ as 1-forms on the orbit space. The bilinear
form degenerates on the discriminant $\Sigma\subset M$ consisting of all
nonregular orbits. 

The Euler vector field is defined by
$$
E={1\over h} \sum_{i=1}^n x_i {\partial\over \partial x_i}
\eqno(2.43)
$$
where $x_1$, \dots, $x_n$ are Euclidean coordinates in ${\bf R}^n$ and $h$
is the Coxeter number of $W$. Recall [Bour] that $h$ is the maximum of the
degrees of basic invariant polynomials $p_1(x)$, \dots, $p_n(x)$, i.e.,
such homogeneous polynomials that
$$
\C [x_1, \dots, x_n]^W = \C [p_1, \dots, p_n]
$$
(Chevalley theorem). Let $\deg \, p_1(x)=h$. Introduce a vector field on
the orbit space
$$
e={\partial\over \partial p_1}.
\eqno(2.44)
$$
It is well-defined up to multiplication by a nonzero constant factor. It
was proved by K.Saito [Sa] (see also [SYS]) that the metric
$$
(~,~)_2 :=\L_e (~,~)_1
\eqno(2.45)
$$
is flat and it does not degenerate globally on $M$. The flat coordinates
of this metric give a distinguished system of generators in the ring of
$W$-invariant polynomials on ${\bf R}^n$ first discovered in [SYS]. In
[Du3]
it was shown (see also [Du4]) that the metrics $(~,~)_1$ and $(~,~)_2$
form
a flat
quasihomogeneous regular pencil of the degree 
$$
d=1-{2\over h}.
\eqno(2.46)
$$
The vector fields $E$ and $e$ have the form (2.43), (2.44), the function
$\tau$ is
$$
\tau={1\over 2\, h} (x,x).
\eqno(2.47)
$$
This produces a polynomial Frobenius structure on the orbit space [{\it
ibid.}].

This construction was generalized in [DZ1] to produce a Frobenius
structure
on orbit spaces of certain extensions of affine Weyl groups. In this case
$d=1$ but the arguments of the above Remark work.
\bigskip
\centerline{\bf 3. Flat pencils and bihamiltonian structures on loop
spaces}
\medskip
We define loop space $L(M)$ of all smooth maps
$$
S^1\to M.
$$
In a local coordinate system $x^1$, \dots, $x^n$ any such a map is given
by a $2\pi$-periodic smooth vector-function $(x^1(s), \dots, x^n(s))$. We
will consider {\it local functionals}
$$
I[x] =\ine P\left( x; \dot x, \ddot x, \dots, x^{(m)}\right) \, ds
\eqno(3.1)
$$
$$
x=\left(x^1, \dots, x^n\right), ~\dot x =\left( \dot x^1, \dots, \dot
x^n\right) = {d x\over ds}, ~\dots
$$ 
as ``functions'' on the loop space. Here $P \left( x; \dot x, \ddot x,
\dots, x^{(m)}\right)$ is a polynomial in $\dot x$, $\ddot x$, \dots
with the coefficients smooth functions of $x$. The nonnegative integer $m$
may depend on the functional.

A {\it local translation invariant Poisson bracket} on the loop space is,
by definition, a Lie algebra structure on this space of functionals
$$
\left( I_1, I_2\right) \mapsto \{ I_1, I_2\}
$$
with the Poisson bracket of the following form
$$
\{ I_1, I_2\} =\ine {\delta I_1\over \delta x^i(s)} A^{ij} {\delta
I_2\over \delta x^j(s)} \, ds
\eqno(3.2)
$$
where $A^{ij}$ is a linear differential operator of some finite order $N$
of the form
$$
A^{ij} = \sum_{k=0}^N a^{ij}_k \left( x; \dot x, \ddot x, \dots,
x^{(m_k)}\right) {d^k\over ds^k}.
\eqno(3.3)
$$
The variational derivatives are the functions defined by the usual rule
$$
I[x+\delta x] - I[x] =\ine {\delta I \over \delta x^i(s)} \delta x^i(s) \,
ds + O \left( || \delta x||^2\right).
\eqno(3.4a)
$$
For the local functionals of the form (3.1) the variational derivatives
are obtained by applying Euler - Lagrange operator
$$
{\delta I\over \delta x^i(s)} ={\partial P\over \partial x^i} -{d\over ds} 
{\partial P \over \partial \dot x^i} + {d^2 \over ds^2} {\partial P\over
\partial \ddot x^i}
- \dots .
\eqno(3.4b)
$$
The coefficients $a^{ij}_k \left( x; \dot x, \ddot x, \dots,
x^{(m_k)}\right)$ of the Poisson bracket must be polynomials in the
derivatives $\dot x$, $\ddot x$ \dots with smooth in $x$ coefficients.

In computations with local Poisson brackets it is convenient to use the
formalism of $\delta$-functions introducing the matrix of distributions
$$
\{ x^i (s_1), x^j(s_2)\} = \sum_{k=0}^N a^{ij}_k \left( x(s_1); \dot
x(s_1), \ddot
x(s_1), \dots, x^{(m_k)}(s_1)\right) \delta^{(k)}(s_1-s_2).
\eqno(3.5)
$$
Here $\delta(s)$ is the $\delta$-function on the circle defined by the
identity
$$
\ine f(s) \delta(s) \, ds = f(0)
$$
for any smooth $2\pi$-periodic function $f(s)$. The derivatives of
$\delta$-function are defined in the standard way by the equations
$$
\ine f(s_2) \delta^{(k)} (s_1-s_2) \, ds_2 = f^{(k)}(s_1).
$$
The formula (3.2) for the Poisson bracket can be recasted into the form
$$
\{ I_1, I_2\} ={1\over 4 \pi^2} \int_0^{2\pi} \int_0^{2\pi} {\delta
I_1\over \delta x^i(s_1)} \{ x^i (s_1), x^j(s_2)\} {\delta I_2 \over
\delta x^j (s_2)} \, ds_1\, ds_2.
\eqno(3.6)
$$

To make the Poisson bracket independent on the choice of the local
coordinates on $M$ the coefficients $a^{ij}_k$ must transform in an
appropriate way with the changes of coordinates $y=y(x)$. The
transformation law of the coefficients is determined by the Leibnitz
identity for the Poisson bracket (3.5)
$$
\{ y^p(s_1), y^q(s_2)\} = {\partial y^p\over \partial x^i}(s_1) {\partial
y^q\over \partial x^j}(s_2) \{ x^i (s_1), x^j(s_2)\}
\eqno(3.7)
$$
together with the following identities for the derivatives of
$\delta$-function
$$
f(s_2) \delta^{(k)}(s_1-s_2) = \sum_{l=0}^k {k\choose l} f^{(l)} (s_1)
\delta^{(k-l)}
(s_1-s_2).
\eqno(3.8)
$$

The constraints for the coefficients $a^{ij}_k$ imposed by skew symmetry
and by Jacobi identity for the Poisson bracket (3.2) can be written as a
finite system of equations of degree one and two resp. for these
coefficients and their derivatives. 

Let us assign degrees to the derivatives putting
$$
\deg {d^k x^i\over ds^k} =k, ~k=1, \, 2, \, \dots
\eqno(3.9a)
$$
and we put
$$
\deg\, f(u)=0
\eqno(3.9b)
$$
for any function independent on the derivatives.
\smallskip
{\bf Definition 3.1.} We say that the bilinear operation (3.2) (or (3.5))
is {\it graded homogeneous of the degree $D$} if the coefficients 
are graded homogeneous polynomials in the derivatives of the degrees
$$
\deg \ a^{ij}_k \left( x; \dot x, \ddot x, \dots, x^{(m)}\right) =D-k, ~
k=0, \, 1, \, \dots
\eqno(3.10)
$$
\medskip
Clearly the order $N$ of (3.5) cannot be greater than the degree $D$.
\smallskip
{\bf Lemma 3.1.} {\it The degree $D$ does not depend on the choice of
local coordinates $x^1$, \dots, $x^n$.}

Proof easily follows from the transformation property (3.7) together with
(3.8).
\medskip
{\bf Example 3.1.} The graded homogeneous Poisson bracket of degree 0 has
the form
$$
\{ x^i(s_1), x^j(s_2)\} = h^{ij}(x) \delta (s_1-s_2)
\eqno(3.11)
$$
where $h^{ij}(x)$ is a usual (i.e., a finite-dimensional one) Poisson
bracket on the manifold $M$.
\smallskip
{\bf Example 3.2.} The graded homogeneous Poisson bracket of degree 1 has
the form
$$
\{ x^i (s_1), x^j(s_2)\} = g^{ij}\left( x(s_1)\right) \dot \delta(s_1-s_2)
+ \Gamma_k^{ij}(x) \dot x^k \delta (s_1-s_2).
\eqno(3.12)
$$
The coefficients $g^{ij}(x)$ and $\Gamma^{ij}_k(x)$ are some functions on
$M$ depending on the choice of local coordinates. This class of Poisson
brackets is called {\it Poisson brackets of hydrodynamics type}. It
was first introduced and studied in [DN1]. The following main result was
proved in this paper (see also [DN3]).

Let us assume that the matrix $g^{ij}(x)$ does not degenerate on an open
dense subset of $M$. (This assumption does not depend on the choice of
local coordinates on $M$.) In this case we call the bracket (3.12) {\it
nondegenerate}.
\smallskip
{\bf Theorem 3.1.} {\it The graded homogeneous nondegenerate Poisson
brackets of the degree 1 on the loop space $L(M)$ are in 1-to-1
correspondence with flat contravariant metrics $g^{ij}(x)$ on $M$. The
coefficients $\Gamma^{ij}_k(x)$ in (3.12) must be the Levi-Civita
contravariant connection of this metric.}
\medskip
{\bf Remark 3.1.} The flat coordinates $t^1$, \dots, $t^n$ of the flat
metric give the densities of Casimirs of the Poisson bracket (3.12)
$$
C^\alpha = \ine t^\alpha(s) \, ds, ~\alpha = 1, \, \dots, n
\eqno(3.13a)
$$
$$
\{ I, C^\alpha \} = 0 ~{\rm for ~ any ~ functional} ~I.
\eqno(3.13b)
$$
Doing a change of the dependent variables
$$
x^i(s) \mapsto t^\alpha(x(s))
$$
we rewrite the Poisson bracket (3.12) in the following constant form
$$
\{ t^\alpha(s_1) , t^\beta(s_2)\} =\eta^{\alpha\beta} \dot \delta(s_1-s_2)
\eqno(3.14)
$$
where the constant coefficients $\eta^{\alpha\beta}$ are the entries of
the matrix of the metric in the flat coordinates $t$.
\medskip
The coefficients of graded homogeneous Poisson brackets of the degree
$D>1$ are also certain differential-geometric objects on the manifold $M$.
These Poisson brackets were first introduced in [DN2] under the name {\it
homogeneous differential-geometric Poisson brackets} (see also [DN3]).

We recall now the general definition of a compatible pair of Poisson
brackets [Mag].
\smallskip
{\bf Definition 3.2.} Two Poisson brackets $\{~,~\}_1$ and $\{~,~\}_2$ on
the same space are called {\it compatible} if the linear combination
$$
\{~,~\}_1-\lambda \{~,~\}_2
\eqno(3.15)
$$
is a Poisson bracket for any $\lambda$.
\medskip
Given a compatible pair of Poisson brackets one can construct certain
family of commuting Hamiltonians. The Hamiltonians $H^{\alpha, p}$ are
determined by the recursion relations
$$
\{ \, . \, , H^{\alpha,p}\}_1 =\{ \, . \, , H^{\alpha,p+1}\}_2, ~~p=0, \,
1, \, \dots, ~\alpha=1, \, \dots, n
\eqno(3.16a)
$$
starting from the Casimirs
$$
H^{\alpha,0} = C^\alpha, ~\alpha=1, \dots, n
\eqno(3.16b)
$$
of the second Posson bracket. By the construction the correspondent
evolutionary systems
admit a bi-hamiltonian structure
$$
{\partial x^i\over \partial T^{\alpha, p}} =\{ x^i(s), H^{\alpha, p}\}_2 =
\{ x^i(s), H^{\alpha, p-1}\}_1, ~~\alpha=1, \, \dots, n, ~~p=1, \, 2,\,
\dots .
$$
In some cases it is possible to prove complete integrability of the
bi-hamiltonian systems.

We prove now the following simple
\smallskip
{\bf Theorem 3.2.} {\it Two graded homogeneous nondegenerate Poisson
brackets of the degree 1 on the loop space $L(M)$ are compatible {\rm iff}
the correspondent flat metrics form a flat pencil.}

Proof. The linear combination (3.15) of two Poisson brackets of the form
(3.12) reads
$$
\eqalign{
& \{ x^i(s_1), x^j(s_2) \}_1 - \lambda \{ x^i(s_1), x^j(s_2)\}_2  \cr
&=\left[ g^{ij}_1(x(s_1)) -\lambda g_2^{ij}(x(s_1))\right]
\dot \delta(s_1-s_2) +
\left[ {\Gamma_1}_k^{ij}(x) -\lambda{\Gamma_2}_k^{ij}(x)\right] \, \dot
x^k \delta (s_1-s_2).
\cr}
\eqno(3.17)
$$
Now the proof immediately follows from Theorem 3.1.
\medskip
{\bf Corollary 3.1.} {\it The loop space $L(M)$ of any Frobenius manifold
$M$ carries a graded homogeneous of degree 1 nondegenerate bi-hamiltonian
structure.}

This follows from Theorem 1.1.
\medskip
Observe that, for $d\neq 1$, the variable
$$
T(s) : =  {2\over 1-d} \tau (s)
\eqno(3.18)
$$
where the flat coordinate $\tau$ was defined in (1.6) has the Poisson
bracket with itself of the form
$$
\{ T(s_1), T(s_2)\}_1 = \left[ T(s_1) + T(s_2)\right] \dot
\delta(s_1-s_2).
\eqno(3.19)
$$
This coincides with the Poisson bracket on the dual space to the Lie
algebra of one-dimensional vector fields (i.e., the Virasoro algebra with
zero central charge). Other Poisson brackets of $T(s)$ are of the form,
due to (2.21), (2.23)
$$
\{ t^\alpha(s_1), T(s_2)\}_1 ={2\over 1-d} E^\alpha(t(s_1)) \dot \delta
(s_1-s_2) + \dot t^\alpha  \delta (s_1-s_2).
\eqno(3.20)
$$
Recall that $E(t)$  depends linearly on $t$.

From the results of Section 2 above it follows that, under the assumption
of quasihomogeneity and regularity, bihamiltonian structures (3.17) on
the loop space $L(M)$ are in
1-to-1 correspondence with Frobenius structures on $M$.

The role of the quasihomogeneity condition in the theory of the degree 1
bihamiltonian structures on $L(M)$ could seem more motivated from the
point of view of a general differential-geometric approach to classical
$W$-algebras outlined in [DZ2]. In this approach we consider Poisson
brackets of the form of formal series in an independent variable
$\epsilon$
$$
\{ x^i(s_1), x^j(s_2)\} =\sum_{k\geq 0} \epsilon^{2k}
\{ x^i(s_1),x^j(s_2)\}^{(k)}
\eqno(3.21)
$$
where the $k$-th coefficient $\{~,~\}^{(k)}$ must be a graded homogeneous
operation of the degree $2k+1$. The skew symmetry and Jacobi identity for
the bracket (3.21) must fulfill as an identity for formal power series
in $\epsilon^2$. The main requirement is that the Poisson bracket (3.21)
must be reducible to the constant form (3.14) by a transformation
$$
x^i =x^i(t) + \sum_{k\geq 1} \epsilon^k Q^i_k (t; \dot t,
\ddot t,
\dots )
\eqno(3.22)
$$
where the coefficients $ Q^i_k (t; \dot t,
\ddot t,
\dots )$ must be graded homogeneous polynomials of the degree $k$ in the
derivatives $\dot t$, $\ddot t$, \dots. Particularly, the leading term
$\{ ~,~\}^{(0)}$ of (3.21) is a graded homogeneous Poisson bracket of
degree 1. Validity of Jacobi identity for (3.21) implies that $\{
~,~\}^{(0)}$ is a Poisson bracket. So, under the nondegeneracy condition
for this bracket, the leading term in (3.22) is given by the flat
coordinates of the correspondent contravariant metric.

We also bring attention of the reader to the construction of [DN1]
(justified in a recent paper [Mal]) of ``averaged Poisson brackets" used
to
describe Hamiltonian structure of Whitham equations. Particularly,
according to this construction, the leading term in the small dispersion
expansion of an arbitrary local Poisson bracket posessing of a
sufficiently rich family of commuting local Hamiltonians is always
given by a degree 1 graded homogeneous Poisson bracket. So, 
(3.21) can be considered as the full small dispersion expansion of the
original Poisson bracket.

$W$-algebras were discovered by A.Zamolodchikov [Za] in order to describe
additional symmetries of conformal field theories with spin greater than
$1/2$. It was realized by Fateev and Lukyanov [FL] that the semiclassical
limit of $W$-algebras coincides with the second Poisson bracket of Gelfand
- Dickey integrable hierarchy. These semiclassical limits of $W$-algebras
were constructed for all simple Lie groups using Drinfeld -
Sokolov construction of the corresponding integrable hierachies [DS].
They were called {\it classical $W$-algebras} (see also [DIZ], [Bouw]).
The role
of the first Poisson bracket of the hierarchy looked not to be relevant in
the construction. However, it will be important in our
differential-geometric approach to classical $W$-algebras and their
generalization. Recall that any Poisson bracket (3.21) by the
assumption has no invariants
w.r.t. the transformations of the form (3.22).

By our definition (see [DZ2]) a classical $W$ algebra is a pair of Poisson
brackets of the form (3.21) such that the linear combination
$\{~,~\}_1-\lambda\{~,~\}_2$ for any $\lambda$ is again a Poisson
bracket
satisfying the above reducibility condition. We also require
validity of certain
quasihomogeneity conditions for the coefficients of the Poisson brackets.
We begin with the leading terms $\{~,~\}_1^{(0)}$ and $\{~,~\}_2^{(0)}$.
The requirement is that this compatible pair of the degree 1 Poisson
brackets corresponds to a quasihomogeneous regular pencil of metrics on
$M$. To motivate this requirement we recall that, according to (3.19),
(3.20) the Poisson bracket  $\{~,~\}^{(0)}_1$ is a nonlinear chiral
extension of the conformal Virasoro algebra with the central charge 0.
The nonzero central charge will arrive with the $\epsilon^2$-correction
(see below). 

Using Theorem 2.2 we see that the leading term in the
$\epsilon^2$-expansion of a classical $W$-algebra is determined by a
Frobenius structure on $M$. Let $E$, $e$, $d$ be resp. the Euler and the
unity vector fields and
the charge of the Frobenius manifold. We will write the quasihomogeneity
conditions for the coefficients of the Poisson brackets using the flat
coordinates $t^\alpha$ on the Frobenius manifold.

Let
$$
\{ t^\alpha (s_1), t^\beta(s_2) \}_1 ^{(k)} 
=\sum _l a^{\alpha\beta}_{k,\,l} (t; \dot t, \ddot t , \dots )
\delta^{(l)}(s_1-s_2)
\eqno(3.23)
$$
$$
\{ t^\alpha (s_1), t^\beta(s_2) \}_2 ^{(k)}
=\sum _l b^{\alpha\beta}_{k,\,l} (t; \dot t, \ddot t , \dots )
\delta^{(l)}(s_1-s_2).
\eqno(3.24)
$$

1. We require that
$$
\L_e a^{\alpha\beta}_{k,\,l} = b^{\alpha\beta}_{k,\, l}
\eqno(3.25)
$$
$$
\L_e  b^{\alpha\beta}_{k,\, l}
=0,
\eqno(3.26)
$$

2. Let us introduce the prolungated vector field 
$$
{\cal E} = E - \sum_{m\geq 1} \sum_{\alpha, \, \beta} \left( m \,
\delta_\beta^\alpha + K^\alpha_\beta\right) {t^\beta}^{(m)} {\partial\over
\partial {t^\alpha}^{(m)}}
\eqno(3.27)
$$
(the matrix $K_\beta^\alpha$ was introduced in (2.17)). Then we require
that
$$
\L_{\cal E} a^{\alpha\beta}_{k, \, l} =(k(d-3) + l)
a^{\alpha\beta}_{k,\,l} + \Lambda^\alpha_\epsilon a^{\epsilon\beta}
_{k,\,l} + a^{\alpha\epsilon}_{k,\, l} \Lambda_\epsilon^\beta
\eqno(3.28)
$$
$$
\L_{\cal E} b^{\alpha\beta}_{k, \, l} =(k(d-3) + l-1)
b^{\alpha\beta}_{k,\,l} + \Lambda^\alpha_\epsilon b^{\epsilon\beta}
_{k,\,l} + b^{\alpha\epsilon}_{k,\, l} \Lambda_\epsilon^\beta.
\eqno(3.29)
$$

3. The first Poisson bracket of the field $T(s)$ given by (3.18) has the
Virasoro
form
$$
\{ T(s_1), T(s_2)\}_1 = [T(s_1)+ T(s_2)] \dot\delta (s_1-s_2) +
\epsilon^2 {c\over 12}  \delta^{(3)} (s_1 -s_2) + O(\epsilon^4).
\eqno(3.30)
$$
The number $c$ is called {\it central charge } of the classical
$W$-algebra.
\medskip
For the clasical $W$-algebras corresponding to the simple Lie groups the
sums in (3.21) are finite. All the coefficients $a^{\alpha\beta}_{k,\,l}$
are polynomials also in $t$. (Observe that in our notations the first and
the second Poisson structures are interchanged.) The central charge $c$ is
equal [FL] to
$$
c=12 \rho^2
\eqno(3.31)
$$
where $\rho$ is half of the sum of positive roots of the root system of
the Lie algebra.  

In [DZ2] it was shown that for any semisimple Frobenius manifold there
exists a germ of order 1 (i.e., the first two terms in (3.21)) of a
classical $W$ algebra with the central charge
$$
c={12\over (1-d)^2} \left[ {n\over 2} - 2 {\rm tr}\, \Lambda^2\right].
\eqno(3.32)
$$
Remarkably, for the Frobenius manifolds on the orbit spaces of
simply-laced Weyl groups (see Example 2.1 above) the formulae (3.31) and
(3.32) give the same result! 

The corrections $\{~,~\}_1^{(1)}$ and
$\{~,~\}_2^{(1)}$ are uniquely
determined by the axioms of Dijkgraaf - Witten [DW] and of Getzler [Ge]
from two-dimensional topological field theory (see explicit formulae in
[DZ2]). The structure of higher order corrections to these brackets
remains unknown. Understanding of this structure could clarify the
eventual role of Frobenius manifolds in the problem of classification of
integrable hierarchies. It will also solve the problem of the genus
expansion in topological field theories (see discussion of this problem in
[DZ2]).
\bigskip
{\bf Acknowledgments.} I wish to thank the organizers of 1997 Taniguchi
Symposium ``Integrable Systems and Algebraic Geometry'' for creative
environment during the sessions and for a generous support.
\vfill\eject 
\centerline{\bf References} \medskip
\item{[Arn]} V.I.Arnol'd, Wave front evolution and equivariant Morse
lemma,
{\sl Comm. Pure Appl. Math.} {\bf 29} (1976) 557 - 582.
\medskip
\item{[Bour]} N.Bourbaki, Groupes et Alg\`ebres de Lie, Chapitres 4, 5 et
6,
Masson, Paris-New York-Barcelone-Milan-Mexico-Rio de Janeiro, 1981.
\medskip
\item{[Bouw]} P.Bouwknegt, K.Schoutens, W-Symmetry, Singapore, World
Scientific, 1995.
\medskip
\item{[DIZ]} P.Di Francesco, C.Itzykson, J.-B.Zuber, Classical
$W$-algebras,
{\sl Comm. Math. Phys.} {\bf 140} (1991) 543 - 567.
\medskip
\item{[DW]} R.Dijkgraaf,
E.Witten, Mean field theory,
topological field theory, and multi-matrix models, {\sl Nucl.
Phys.} {\bf B342} (1990), 486--522.
\medskip
\item{[DS]} V.Drinfeld, V.Sokolov, Lie algebras and equations of Korteweg
- de Vries type, {\sl J. Sov. Math.} {\bf 30} (1985) 1975 - 2036.
\medskip
\item{[Du1]} B.Dubrovin, Integrable systems in topological
field theory, {\sl Nucl. Phys.} {\bf B379} (1992), 627--689.
\medskip
\item{[Du2]} B.Dubrovin, Topological conformal
field theory from the point of view of integrable systems,
In: Integrable Quantum Field Theories, Edited
by L.Bonora, G.Mussardo, A.Schwimmer, L.Girardello, and M.Martellini,
Plenum Press, NATO ASI series {\bf B310} (1993) 283 - 302.
\medskip
\item{[Du3]} B.Dubrovin, Differential geometry of the space of orbits
of a Coxeter group, Preprint SISSA-29/93/FM (February 1993),
hep-th/9303152.
\medskip
\item{[Du4]} B.Dubrovin, Geometry of 2D topological field theories,
in: Integrable Systems and Quantum Groups, Montecatini Terme, 1993.
Editor: M.Francaviglia, S. Greco. Sprin\-ger Lecture Notes in
Math. {\bf 1620} (1996), 120--348.
\medskip
\item{[DFN]} B.Dubrovin, A.T.Fomenko, S.P.Novikov, Modern Geometry - 
Methods and Applications. v.1, New York, Springer-Verlag, 1984.
\medskip
\item{[DN1]} B.Dubrovin and S.P.Novikov, The Hamiltonian formalism of
one-dimensional systems of the hydrodynamic type and the Bogoliubov -
Whitham averaging method, {\sl Sov. Math. Dokl.} {\bf 27} (1983) 665 -
669.
\medskip
\item{[DN2]} B.Dubrovin and S.P.Novikov, 
On Poisson brackets of hydrodynamic type,
{\sl Soviet Math. Dokl.} {\bf 279} (1984), 294-297. 
\medskip
\item{[DN3]} B.Dubrovin and S.P.Novikov, Hydrodynamics of weakly deformed
soliton lattices. Differential
geometry and Hamiltonian theory,
{\sl Russian Math. Surv.} {\bf 44:6} (1989), 29-98.
\medskip
\item{[DZ1]} B.Dubrovin and Youjin Zhang, Extended affine Weyl groups and
Frobenius manifolds, {\sl Compositio Math.}, {\bf 111} (1998) 167-219.
\medskip
\item{[DZ2]} B.Dubrovin and Youjin Zhang, Bihamiltonian hierarchies in
2D 
topological field theory at one-loop approximation, Preprint
SISSA 152/97/FM, hep-th/9712232.
\medskip
\item{[FL]} V.A.Fateev, S.L.Lukyanov,
Additional symmetries
and exactly-solvable models in two-dimensional
conformal field theory, parts I, II, and III,  
{\sl Sov. Sci. Rev.} {\bf A15} (1990) 1.
\medskip
\item{[Ge]} E.Getzler, Intersection theory on ${\bar M}_{1,4}$
and elliptic Gromov-Witten invariants, alg-geom/9612004, to appear
in J. Amer. Math. Soc..
\medskip
\item{[Mag]} F. Magri, A simple model of the integrable Hamiltonian
systems,
{\sl J. Math. Phys.} {\bf 19}(1978), 1156-1162.
\medskip
\item{[Mal]}  A.Ya.Maltsev,  The conservation of the Hamiltonian
structures in Whitham's
method of averaging, Preprint, solv-int/9611008. 
\medskip
\item{[Sa]} K.Saito, On a linear structure of a quotient variety
by a finite
reflection group, Preprint RIMS-288 (1979).
\medskip
\item{[SYS]} K.Saito, T.Yano, J.Sekeguchi, On a certain generator
system
of the ring of invariants of a finite reflection group,
{\sl Comm. in Algebra} {\bf 8(4)} (1980) 373 - 408.
\medskip
\item{[Za]} A.Zamolodchikov, Additional infinite symmetries in
two-dimensional conformal quantum field theory, {\sl Theor. Math. Phys.}
{\bf 65} (1985) 1205 - 1213. 
\vfill\eject\end